\documentclass[12pt,a4paper]{amsart}
\usepackage{amsmath}\usepackage{amssymb}
\usepackage{amsthm}
\usepackage[english]{babel}\usepackage{german}
\usepackage{nicefrac}  \usepackage{floatflt}
\usepackage{graphicx}    
\usepackage{xcolor}
\graphicspath{ {./Images/} }
\usepackage{xr}

\usepackage{biblatex}
\author{Lauritz Streck} \thanks{The author has received funding from the European Research Council (ERC)
under the European Union’s Horizon 2020 research and innovation programme (grant agreement No. 803711)}
\date{\today}
\title[On equidistribution of polynomial sequences]{On equidistribution of polynomial sequences in quotients of $PSL_2(\mathbb{R})$}

\begin{document}
\newtheorem{theorem}{Theorem}
\newtheorem*{theorem*}{Theorem}
\newtheorem{definition}[theorem]{Definition}
\newtheorem{remark}[theorem]{Remark}
\newtheorem{lemma}[theorem]{Lemma}
\newtheorem{corollary}[theorem]{Corollary}
\newtheorem{prop}[theorem]{Proposition}
\newtheorem{claim}[theorem]{Claim}
\newtheorem{observation}[theorem]{Observation}
\newcommand{\eps}{\varepsilon}
\newcommand{\R}{\mathbb{R}}
\newcommand{\Z}{\mathbb{Z}}
\newcommand{\N}{\mathbb{N}}
\newcommand{\Imag}{\mathrm{Im}}
\newcommand{\abcd}{\begin{pmatrix}
a & b \\ c & d
\end{pmatrix}}

\begin{abstract}
In this paper,  it is shown that for every lattice $\Gamma \subset  PSL_2(\mathbb{R})$ there exists a $c>0$ such that for any $0 \leq \gamma<c$ the sequence $p h(n^{1+\gamma})$ equidistributes for any $p  \in \Gamma \backslash PSL_2(\mathbb{R})$, where $h$ is the horocycle flow. This makes modest progress towards a conjecture of Shah and generalizes a result of Venkatesh („Sparse equidistribution problems, period bounds, and subconvexity“, 2005),  who established the same equidistribution for co-compact lattices.  The proof utilizes a dichotomy between good equidistribution estimates and approximability of $\{p h(t), t \leq T \}$ by closed horocycles of small period.  
\end{abstract}
\maketitle
\section{Introduction}
Consider the (multiplicative) group $G:=PSL_2(\R)$ with a Haar measure $\mu_G$. A \textit{lattice} $\Gamma \subset G$ is a discrete subgroup such that the quotient $X:=\Gamma \backslash G$ has a fundamental domain in $G$ of finite Haar measure.  The Haar measure then descends to a finite measure $\mu_X$. We define the matrices 
\[
h(x):=\begin{pmatrix}
1 & x \\ 0 & 1
\end{pmatrix} \quad
a(y):=\begin{pmatrix}
y^{\frac{1}{2}} & 0 \\ 0 & y^{-\frac{1}{2}}
\end{pmatrix}.
\]
The \textit{geodesic flow} at time $t$ of $p \in X$ is defined by $g_t(p):=p  a(e^t)$ and the \textit{horocycle flow} at time $t$ is defined by $h_t(p ):=p  h(t)$.  

While the orbit $g_t(p )$ for $t \to \infty$ can behave quite irregularly depending on the initial point, the horocycle orbit $h_t(p )$ is known to behave much more rigidly.  Before we detail the known results, we pin down some notation. We say that the orbit $h_t(p )$ equidistributes with respect to $\mu_X$ if for any compactly supported, continuous function $f$ on $X$, 
\[
\lim_{T \to \infty} \frac{1}{T} \int_{0}^T f(ph(t)) \; dt \to \int f \; d\mu_X.
\]
Similarly, we say that the orbit equidistributes along a sequence $a_n \in \R$ with respect to $\mu_X$ if 
\[
\lim_{N \to \infty} \frac{1}{N} \sum_{n=0}^{N-1} f(ph(a_n)) \; dt \to \int f \; d\mu_X.
\]
Lastly, a point $p  \in X$ is called \textit{periodic} if there is a $t_0 \in \R$ such that $p =p h(t_0)$.  In this case, the horocycle orbit will be trapped in the periodic orbit and will never equidistribute with respect to $\mu_X$; the system $t \mapsto p h(t)$ is then isomorphic to the circle-rotation $x \mapsto x+t_0^{-1}$ on the torus $\R/\Z$.  Below, we use ``$p h(a_n)$ equidistributes'' as a shorthand for ``for all non-periodic $p \in X$,  $p h(a_n)$ equidistributes with respect to $\mu_X$''.   It was shown by Dani and Smillie that both $p  h(t)$ for $t \in \R$ and $p  h(n)$ for $n \in \mathbb{N}$ equidistribute.  

It was subsequently asked what happens for sequences other than $\N$.  Margulis conjectured that $p h(p_n)$, where $p_n$ is the $n$-th prime number,  should also equidistribute.  Shah conjectured that for any $\gamma \geq 0$,  $p h(n^{1+\gamma})$ would equidistribute.  We remark that these results follow for $\mu_X$-almost every $p \in X$ from the work of Bourgain in a much more general context \cite{bourgain}.  The challenge is really to establish equidistribution for \textit{all} non-periodic $p \in X$.

Venkatesh made progress on Shah's conjecture by showing that for co-compact $\Gamma$,  there is a small $c=c(\Gamma)>0$ such that for all $0 \leq \gamma<c$ and all $p \in X$,  $p h(n^{1+\gamma})$ equidistributes \cite{venkatesh}.  His proof operates by controlling arithmetic sequences of the type $p h(sn)$ for $n \in \{0, \dots, N-1\}$ with $s$ small compared to $n$.  Controlling these sparse sequences also means that the almost-primes equidistribute for co-compact $\Gamma$; that is,  for sufficiently big $R$,  $p h(q)$ equidistributes,  where $q$ runs over all numbers having at most $R$ many prime factors.  That controlling sparse sequences is enough to control the almost-primes can be seen either using sieve methods or using the pseudo-random measure $\nu$,  introduced by Goldston and Yilmaz and subsequently used by Green and Tao to show that the primes contain infinitely long arithmetic progressions \cite{goldston}, \cite{green-tao} (see \cite{sarnak-ubis} for a proof of the equidistribution of almost-primes using sieve methods and \cite{streck} for a proof using the pseudo-random measure $\nu$).  

Sarnak and Ubis showed that the almost-primes equidistribute for $\Gamma=PSL_2(\Z)$, which is not co-compact \cite{sarnak-ubis}.  It was subsequently proved by the author that the almost-primes equidistribute for all lattices $\Gamma$ in $PSL_2(\R)$ \cite{streck}.  

In this paper,  the equidistribution of $p h(n^{1+\gamma})$ is established for small $\gamma$ in the setting of a general lattice. This generalises Venkatesh's result from co-compact $\Gamma$ to all lattices $\Gamma$ in $PSL_2(\R)$ and makes modest progress on the conjecture of Shah. 

We make this precise in the result below, which is the main result of this paper. For this, we need some more notation and start by defining the metric $d_X$. The group $G=PSL_2(\R)$ comes with a natural left-invariant metric $d_G$ (see for example Chapter 9 in \cite{einsiedler-ward}). This metric descends to $X$ via $d_X(\Gamma g, \Gamma h):=\inf_{\gamma \in \Gamma} d_G(g, \gamma h)$.  We also fix a point $p_0 \in X$ and define $\mathrm{dist}(p):=d_X(p, p_0)$.  

For two functions $f, g \colon U \to \R$,  we write $f \ll g$ or $f=O(g)$ if there is a constant $C$ such that $|f(x)| \leq C |g(x)|$ for all $x \in U$, where $U$ is some domain.   In this paper,  this constant $C$ implicit in the definition is always allowed to depend on the lattice $\Gamma$ and the choice of $\gamma$,  but nothing else.  We write $f \sim g$ if both $f \ll g$ and $g \ll f$.

For a function $f \in C^4(X)$,  let $\Vert f \Vert_{W^4}$ be its Sobolev norm in the Hilbert space $W^{4, 2}$ involving the fourth derivative, and let $\Vert f \Vert_{\infty, j}$ be the supremum norm of the $j$-th derivatives. Define
\[
\Vert f \Vert:=\Vert f \Vert_{W^4}+\Vert f \Vert_{\infty, 1}+\Vert f \Vert_{\infty, 0};
\]
this norm is the same one Str\"ombergsson used to show his equidistribution result \cite{strombergsson}.  We let $\beta$ be the constant in Theorem \ref{thm::venkatesh}; it ultimately comes from the rate of effective mixing.  The constant in Theorem \ref{powerequidist} can be taken to be $c=\frac{\beta}{600}$.
\begin{theorem} \label{powerequidist}
For any lattice $\Gamma \subset PSL_2(\mathbb{R})$ there is a constant $c=c(\Gamma)>0$ such that for any $0 \leq \gamma \leq c$,  any non-periodic $p \in X$ and any function $f \in C^4(X)$ with $\Vert f \Vert=1$,
\[
\left| \frac{1}{T} \sum_{n \leq T} f\left(p h\left(n^{1+\gamma}\right)\right) -\int f \; d\mu_X \right| \ll r^{-\frac{\beta}{4}},
\]
where $r=T^{1+\gamma} \exp(-\mathrm{dist}(g_{\log T^{1+\gamma}}(p)))$. Because $r \to \infty$ as $T \to \infty$, the sequence $p h(n^{1+\gamma})$ equidistributes.
\end{theorem}

To prove Theorem \ref{powerequidist},  we will split the range into different intervals and use Taylor expansion on each one. On an interval $[T_0,  T_1]$,  the function $t^{1+\gamma}$ will be approximately equal to $T_0^{1+\gamma} + (1+\gamma) T_0^\gamma (t-T_0)$, provided that $T_0$ is not too small and that the range is not too long.  The question thus becomes how well $ph(ns)$ for $s \sim T^{\gamma}$ equidistributes.  To control these sparse arithmetic sequences,  we need two results.

The first one is the following theorem, which is a straightforward consequence of combining Str\"ombergsson's equidistribution result \cite{strombergsson} with Venkatesh's method \cite{venkatesh}, as performed for example by Zheng \cite{zheng}. 
\begin{theorem}[\cite{zheng}, Theorem 1.2] \label{thm::venkatesh}
Let $\Gamma$ be a non-compact lattice in $G$. Let $f \in C^4(X)$ with $\Vert f \Vert < \infty$ and $1 \leq s<T$. Then
\[
\left| \frac{s}{T} \sum_{1 \leq j \leq \nicefrac{T}{s}} f(p h(sj))-\int f \; d\mu_X \right| \ll s^{\frac{1}{2}} r^{-\frac{\beta}{2}} \Vert f \Vert
\]
for any initial point $p \in X$, where $r=T \exp(-\mathrm{dist}(g_{\log T}(p)))$. The parameter $\frac{1}{6}>\beta>0$ and the implied constant depend only on $\Gamma$.
\end{theorem}
In the cases that $r$ is big compared to $T$ (say $r \geq T^\eps$ for some absolute $\epsilon$), this result in itself is enough to show equidistribution of the sequence $ph(n^{1+\gamma})$.  

The result below will be used to deal with the case in which the equidistribution is bad.  It was proved by the author in \cite{streck} in order to show equidistribution of almost-primes.  Its proof uses ideas of Sarnak and Ubis \cite{sarnak-ubis} and has parallels to \cite{strombergsson},  whose proof in turn uses ideas going back to Marina Ratner.  This result encompasses the dichotomy mentioned in the abstract.
\begin{lemma}[Lemma 1.3 in \cite{streck}] \label{lem::orbitapprox}
Let $\Gamma$ be a lattice in $G=PSL_2(\R)$ and let $X=\Gamma \backslash G$.  Let $p \in X$ and $T \geq 0$. Let $\delta>0$ and $K \leq T$. 

There is an interval $I_0 \subset [0,T]$ of size $|I_0| \leq \delta^{-1} K^2$ such that:\\
For all $t_0 \in [0,T] \backslash I_0$,  there is a segment $\{p h(t),  t \leq K\}$ of a closed horocycle approximating $\{ph(t_0+t), 0 \leq t \leq K\}$ of order $\delta$,  in the sense that
\[
\forall 0 \leq t \leq K: \quad d_X\left(ph(t_0+t),  p h(t)\right) \leq \delta. 
\]
The period $P=P(t_0, p)$ of this closed horocycle is at most $ P \ll  r$, where $r=T \exp(-\mathrm{dist}(g_{\log T}(p)))$.

Moreover,  one can assure $P \gg \eta^2 r$ for some $\eta>0$ by weakening the bound on $I_0$ to $|I_0| \leq \max\left(\delta^{-1} K^2,  \eta T\right)$.
\end{lemma}
\subsection*{Acknowledgements}
The following paper is a follow-up paper to \cite{streck}, which is based on the master's thesis I did at the Hebrew University of Jerusalem in 2020.  As such, I am thankful for the support by my thesis advisor Tamar Ziegler and by Elon Lindenstrauss, who also suggested that the result in the present paper should be achievable with the ideas in \cite{streck}.  I thank my PhD supervisor P\'eter Varj\'u for giving me the freedom to finish the work on these two papers while doing my PhD with him. 
Above all, I am grateful to Adri\'an Ubis,  who suggested the argument used in the proof of Claim \ref{claim::qbig} in his review of the previous paper,  simplifying the proof in \cite{streck} considerably.  Without getting this new perspective on the material two years later,  I would not even have thought of revisiting the problem solved in this paper. 

\section{On the behaviour of the equidistribution parameter in Theorem \ref{thm::venkatesh}} \label{sec2}
Except for Lemma \ref{lem::orbitapprox} itself,  we will also need some of the other material in Chapter 4 of \cite{streck} in order to prove Theorem \ref{powerequidist}. We recall some of the material, going slightly beyond what is presented in \cite{streck}.

It is well known that $G \cong T_1 \mathbb{H}$, where $\mathbb{H}$ is the upper half-plane with the hyperbolic metric. $X=\Gamma \backslash G$ then has as fundamental domain a set $T_1F$,  where $F$ is a geodesic polygon in $\mathbb{H}$ - that is, a polygon with finitely many vertices with the edges being pieces of geodesics \cite{einsiedler-ward}.  This fundamental polygon $F$ has finitely many vertices touching the boundary of the upper half-plane, either at the axis with real part equal to zero or at infinity.  After identifying vertices that are in the same orbit under the action of $\Gamma$,  one gets the cusps of $X$, which we will denote by $r_1, \dots,  r_n$.  Any such cusp $r_i$ is in 1-1 correspondence to an element $\gamma_i \in \Gamma$ with the property that $\gamma_i$ fixes $r_i$ and that $\gamma_i$ is conjugated to $h(1)$ (see Lemma 3.1 in \cite{streck}). For each cusp, there are elements $\sigma_i \in G$ such that $\sigma_i r_i=\infty$ and $\sigma_i \gamma_i \sigma_i^{-1}=h(1)$.  

For $g \in G$, we define $Y^0_i(g):=\mathrm{Im}(\sigma_i g)$, where 
\[
\Imag\left(\abcd \right):=\frac{1}{c^2+d^2} 
\]
is the imaginary part of the the matrix projected to $\mathbb{H}$.  We also set for $p=\Gamma g_p \in X$,  $y_i^0(p):=\max_{\gamma \in \Gamma} Y^0_i(\gamma g_p)$. 

It was shown in Lemma 4.1 in \cite{streck} that there exist disjoint neighbourhoods $C_i \subset X$ of each cusp $r_i$ with $K=X \backslash \bigcup C_i$ being compact such that for any $p \in C_i$,  $\exp(\mathrm{dist}(p)) \sim y_i^0(p)$ (while of course $\exp(\mathrm{dist}(p)) \sim 1$ for $p \in K$).  Arguing as in the proof of 1.  in Lemma 4.1,  one also sees that if $p=\Gamma g_p \in C_i$ and $g_p$ is such that $Y_i^0(g_p)=y_i^0(p)$,  then for any $\gamma \in \Gamma$,  either $Y_i^0(\gamma g_p) \ll 1$ or $Y_i^0(\gamma g_p)=Y_i^0(g_p)$ (which is the case in which $\sigma_i \gamma g_p=h(n) \sigma_i g_p$ and $\gamma=(\gamma_i)^n$ for some $n$).  This implies in particular that there is an absolute constant $C=C(\Gamma)$ such that if $g_p$ is such that $Y_i^0(g_p) \geq C$,  then 
\[
Y_i^0(g_p) \sim y_i^0(p) \sim\exp(\mathrm{dist}(p)),
\]
where the second equivalence holds because $y_i^0(p) \geq C$ implies that $p \in C_i$ for $C$ sufficiently big.

We will use the equidistribution parameter $r$ in the statement of Theorem \ref{powerequidist} with varying orbit lengths,  so we set 
\[
r(q, K):=K \exp(-\mathrm{dist}(g_{\log K}(q))).
\]
\begin{observation} \label{obs::coordsystem}
There is an absolute $c_0=c_0(\Gamma)>0$ such that for any $T$ and any $p$,  if there is a representative $g_p$ of $p$ and an $i$
such that for $\sigma_i g_p=:\abcd$,  $\max(T^2 c^2, d^2) \leq c_0 T$,  then $r(p, T) \sim \max(T^2 c^2, d^2)$. 
\end{observation}
\begin{proof}
We have that 
\[
2 \max(T^2 c^2, d^2) \geq T (c^2T+d^2T^{-1})=Y_i^0(g_{\log T}(g_p))^{-1} T.
\]
Thus,  $Y_i^0(g_{\log T}(g_p)) \geq \frac{1}{2} c_0^{-1}$, which shows that 
\[
\exp(\mathrm{dist}(g_{\log T}(p))) \sim Y_i^0(g_{\log T}(g_p))
\]
by the argument above,  provided that $c_0$ is sufficiently small.   
\end{proof}

\section{Proof of Theorem \ref{powerequidist}}
We start by approximating $t^{1+\gamma}$ with sparse arithmetic sequences.  More precisely,  we write 
\[
t^{1+\gamma}=T_0^{1+\gamma}+(1+\gamma) T_0^\gamma (t-T_0)+O(T^{-\frac{1}{6}})
\]
on $[T_0, T_0+T^{\frac{1}{3}}]$ for $T_0 \geq T^{\frac{5}{6}}$ using Taylor expansion. 

We will split into several cases. To govern which case we are in,  we fix some $\varepsilon>0$ and impose that $\gamma < \frac{\eps \beta}{6}$. We will see at the end which value of $\eps$ makes everything work (which will turn out to be $\eps=\frac{1}{100})$.  

To apply the results about sparse equidistribution,  we are thus tasked with evaluating expressions of the form 
\[
\left| \frac{1}{K} \sum_{n \leq K} f\left(q h\left((1+\gamma) T_0^\gamma n \right)\right) -\int f \; d\mu_X \right|
\]
for $q=ph(T_0^{1+\gamma})$ and $T^{\frac{1}{6}} \leq K \leq T^{\frac{1}{3}}$, given some $T_0 \leq T$.  In the case that $r(q, K) \geq T^\eps$,  Theorem \ref{thm::venkatesh} is enough to deduce good equidistribution.

If $r(q, K) \leq T^\eps$,  then $g_{\log K}(q)$ must lie in the neighbourhood $C_i$ of some cusp $r_i$,  as explained in the previous section.  In this case,  there is a (essentially unique) representative $g_q$ of $q$ such that $r(q, K) \sim \max(K^2 c^2, d^2)$,  where we set
\[
\abcd:=\sigma_i g_q,
\]
now and for the next couple of pages.

One then has to split into two more cases. The distinction between these cases is governed by
\[
W_q:=\left| \frac{d}{c} \right|.
\]
The relevance of this $W_q$ is that it measures the time it takes until one gets from bad to good equidistribution again.  More precisely,  by Observation \ref{obs::coordsystem},  
\begin{equation} \label{eq::Wimp}
r(q, K) \sim \begin{cases} d^2,  K \leq W_q \\
d^2 \frac{K^2}{W_q^2},  K \geq W_q
\end{cases}
\end{equation}
as long as $r(q, K) \leq c_0 K$.  

This means that even if $q$ and $K$ are such that $r(q, K) \leq T^{\eps}$,  one has that $r(q, T^\eps W_q) \geq T^{2\eps}$.  Together with Theorem \ref{thm::venkatesh} this will be good enough to show effective equidistribution under all assumptions except for the ones of Proposition \ref{prop::intcase} below.  Under those assumptions,  which encompass the most interesting case,  almost the entire horocycle orbit $\{ph(t),  t \leq T^{1+\gamma}\}$ is close to periodic horocycle orbits of small period.  In this case, one will need Lemma \ref{lem::orbitapprox} to conclude. 

\begin{prop} \label{prop::intcase}
Let $\Gamma$ and $\gamma<c$ be as in Theorem \ref{powerequidist} and let $\eps=\frac{1}{100}$. 
Let $p \in X$ and $T$ be such that $r(p, T^{1+\gamma}) \leq T^{4\eps}$ and $W_p \geq T^{1-\eps}$.  Then for $f$ as in Theorem \ref{powerequidist},
\[
\left|\frac{1}{T} \sum_{n \leq T} f(ph(n^{1+\gamma})) - \int f \; d\mu_X \right| \ll r^{-\frac{\beta}{4}}.
\]
\end{prop}
To prove Theorem \ref{powerequidist}, we will first show how one can reduce its proof to Proposition \ref{prop::intcase} using Observation \ref{obs::coordsystem} and Theorem \ref{thm::venkatesh}. We will then prove Proposition \ref{prop::intcase}.

\begin{proof}[Proof of Theorem \ref{powerequidist} assuming Proposition \ref{prop::intcase}]
Say we are given some $t_0$ and set $q=ph(t_0^{1+\gamma})$.  If $r:=r(q, T^{\frac{1}{6}}) \geq T^{\eps}$, then we know by Theorem \ref{thm::venkatesh} that for any $f$ with $\Vert f \Vert \leq 1$,
\[
\left| \frac{1}{T^{\frac{1}{6}}} \sum_{n \leq T^{\frac{1}{6}}} f\left(q h\left((1+\gamma) t_0^\gamma n \right)\right) -\int f \; d\mu_X \right| \ll T^{\frac{\gamma}{2}} r^{-\frac{\beta}{2}} \leq r^{-\frac{\beta}{4}},
\]
where we recall $\gamma \leq \frac{\eps \beta}{6}$.
We are thus done unless there is a $q$ such that $r=r(q, T^{\frac{1}{6}}) \leq T^{\eps}$.  As we saw in Section \ref{sec2},  then with $c$ and $d$ as defined on the previous page,
\begin{equation} \label{eq:rtilde}
r \sim \max\left(T^{\frac{2}{6}} c^2, d^2\right).
\end{equation}
If $c^2T^{\frac{2}{6}}$ attains the maximum in (\ref{eq:rtilde}), or equivalently, if $W_q \leq T^\frac{1}{6}$,  then $r(q, T^{\frac{1}{4}}) \sim T^{\frac{1}{6}} r \geq T^{\frac{1}{6}}$ by (\ref{eq::Wimp}) and we are done by Theorem \ref{thm::venkatesh}.  We can thus assume $W_q \geq T^{\frac{1}{6}}$.  The claim below shows how one can improve the lower bound on $W_q$ further.
\begin{claim} \label{goodapprox}
Let $q=p h(t_0^{1+\gamma})$ such that $r \leq T^\eps$.  Set $W:=W_q$.  If $W \leq T^{1-\eps}$,  then for $K=W^{1+\eps}$ and for $f$ with $\Vert f \Vert \leq 1$,
\[
\left| \frac{1}{K} \sum_{0 \leq n \leq K} f\left(p h\left((t_0+n)^{1+\gamma}\right)\right) -\int f \; d\mu_X \right| \ll r^{-\frac{\beta}{4}}.
\]
\end{claim}
\begin{proof}[Proof of Claim \ref{goodapprox}]
Fix some $W^{1+\eps} \geq s \geq W^{1+\frac{\eps}{2}}$ and note that then $c^2s^2 \sim W^{-2} s^2 d^2 \gg d^2$. Thus,
\begin{align*} 
r(qh(s), T^{\frac{1}{3}}) &\sim \max \left(T^{\frac{2}{6}} c^2,  (d+cs)^2 \right) \sim \max \left(T^{\frac{2}{6}} c^2,  c^2 s^2 \right)\\
&=c^2 s^2 \sim \left(\frac{s}{W}\right)^2 r \geq r W^\eps \geq r T^{\frac{\eps}{6}},
\end{align*}
where the first equivalence is due to Observation \ref{obs::coordsystem},  which is applicable because $ \left(\frac{s}{W}\right)^2 r \ll T^{3\eps}$.  Applying Theorem \ref{thm::venkatesh} shows that
\begin{align*}
&\left| \frac{1}{T^{\frac{1}{3}}} \sum_{n \leq T^{\frac{1}{3}}} f\left(p h(t_0^{1+\gamma}+s)  h\left((1+\gamma) (t_0+s)^\gamma n \right)\right) -\int f \; d\mu_X \right| \\
&\ll T^{\frac{\gamma}{2}}  T^{-\frac{\epsilon \beta}{12}} r^{-\frac{\beta}{2}}\leq r^{-\frac{\beta}{2}}.
\end{align*}
Now we use Taylor approximation as above to split the orbit of $(t_0+n)^\gamma$ with $n \leq K$ into different ranges $[s, s+T^{\frac{1}{3}}]$ and note that for all but a $W^{-\frac{\eps}{2}} T^{\gamma}$ proportion of $s$, one has $(t_0+s)^{1+\gamma}-t_0^{1+\gamma} \geq W^{1+\frac{\eps}{2}}$.  As
$W^{-\frac{\eps}{2}} T^{\gamma} \leq T^{-\frac{\eps}{12}} \leq r^{-\frac{\beta}{4}}$, the claim is shown.
\end{proof}
We have thus shown the conclusion of Theorem \ref{powerequidist} unless there is a $q=p h(t_0)$ such that $r(q,  T^{\frac{1}{6}}) \leq T^\eps$ and $W_q \geq T^{1-\eps}$.  We let $c$ and $d$ be as defined above and note that in the case considered,  $r(q,  T^{\frac{1}{6}}) \sim \max(c^2 T^{\frac{1}{3}}, d^2)=d^2$ by definition of $W_q$.  By (\ref{eq::Wimp}), this implies that 
\[
r(q, T^{1+\gamma}) \ll  d^2 \frac{T^{2(1+\gamma)}}{W^2_q}\ll T^{4\eps}.
\]
Lastly,  to get an error term in $r(p,  T^{1+\gamma})$ instead of $r(q, T^{1+\gamma})$, we note that $g_q h(-t_0^{1+\gamma})$ is a representative of $p$ and that because $(d-ct_0^{1+\gamma})^2 \ll d^2  T^{2(1+\gamma)} W^{-2} \ll T^{4\eps}$,
\[
r(p, T^{1+\gamma}) \sim \max\left(c^2 T^{2(1+\gamma)},  (d-ct_0)^2 \right) \ll T^ {4 \eps},
\]
where the first equivalence is due to Observation \ref{obs::coordsystem}.  We have thus reduced the proof of Theorem \ref{powerequidist} to the assumptions of Proposition \ref{prop::intcase}.
\end{proof}
It now only remains to show Proposition \ref{prop::intcase}, which is the main part of the proof of Theorem \ref{powerequidist}.  
\begin{proof}[Proof of Proposition \ref{prop::intcase}]
Let $p$ and $T$ be given such $r:=r(p, T^{1+\gamma}) \leq T^{4\eps}$ and $W:=W_p \geq T^{1-\eps}$.  Here,  $W=\left| \frac{d}{c}\right|$,  with $c$ and $d$ as defined in Observation \ref{obs::coordsystem}.  We also let $g:=g_p$ and $\sigma_i$ be as in Observation \ref{obs::coordsystem}. We invoke Lemma \ref{lem::orbitapprox} to split the orbit $[0, T^{1+\gamma}]$ into pieces of length $K=T^{\frac{1}{3}}$.  As in the proof of Lemma 1.3 in Chapter 4 of \cite{streck},  we now parametrize the orbit using the equation 
\[
\sigma_i g h(W+s)=lh(s)=h\left(\alpha-\frac{Rs}{s^2+1}\right)a\left(\frac{R}{s^2+1}\right)k(-\mathrm{arccot} \; s)
\]
where $l:=\sigma_i g h(W)=:(\alpha+iR,  -i)$ is the highest point of the horocycle orbit.  Given an $M \leq T$,  we then have that $p h(M^{1+\gamma}+t), t \leq T^{\frac{1}{3}}$ is at distance at most $O(T^{-\frac{1}{6}})$ from the orbit on a periodic horocycle $\xi h(t), t \leq T^{\frac{1}{3}}$ with its period being equal to $y^{-1}$, where 
\[
y:=\frac{R}{(M^{1+\gamma}-W)^2+1}. 
\]
By the second clause in Lemma \ref{lem::orbitapprox}, we can assume $r \gg y^{-1} \gg \delta^2 r$ except on an interval of proportion $\delta$, where $\delta$ is to be chosen later.
Using Taylor approximation on $t^{1+\gamma}$,  we thus want to bound
\[
\left| \frac{(1+\gamma) M^\gamma}{T^{\frac{1}{3}}}\sum_{(1+\gamma) M^\gamma n \leq T^{\frac{1}{3}}} f(\xi h((1+\gamma)M^\gamma n)) - \int f d\mu_X \right|.
\]
However,  we may run into problems here: If for example $y^{-1}=(1+\gamma)M^\gamma$,  the points do not equidistribute at all in the periodic horocycle.  To deal with this and related obstructions,  we proceed similarly to the proof of Claim 5.2 in \cite{streck}.  For notational convenience,  we set $s:=(1+\gamma) M^\gamma$.  Let $ q \in \mathbb{N}$ with $y^{-1} \leq q \leq y s^{-1} T^{\frac{1}{3}}$ be such that
\[
\left| s y - \frac{a}{q} \right| \leq \frac{y^{-1} s}{q T^{\frac{1}{3}}}
\]
for some $a$ coprime to $q$ (such $q$ exists by the pigeonhole principle).  The problem case occurs if $q$ is small compared to $y^{-1}$.  If on the other hand $q$ is sufficiently big,  there are so many distinct points in the interval $[0, y^{-1}]$ that they cannot help being dense enough to approximate $\int_0^{1} f(\xi h(t y^{-1})) dt$ by force, as we show now.
\begin{claim} \label{claim::qbig}
If $q \geq y^{-3}$,  then
\[
\left| \frac{s}{T^{\frac{1}{3}}}\sum_{s n \leq T^{\frac{1}{3}}} f(\xi h(s n)) -\int_0^{1} f(\xi h(t y^{-1})) dt \right| \ll y \ll \delta^{-2} r^{-1},
\]
where $q,  s,  y$ and $\xi$ all depend on $M$.
\end{claim}
\begin{proof}[Proof of Claim \ref{claim::qbig}]
(The argument in the proof of this claim was suggested by Adri\'an Ubis) We set $F(t):=f(\xi h(t y^{-1}))$, which is one periodic. Because the function $f$ is $1$-Lipschitz with respect to the hyperbolic metric,  the function $F$ is $y^{-1}$-Lipschitz.  We wish to show 
\[
\left| \frac{s}{T^{\frac{1}{3}}} \sum_{sn \leq T^{\frac{1}{3}}} F(nsy)-\int_0^1 F(t) dt \right| \ll y.
\]
For this, we note that as for any $n$
\[
\left|sny-n\frac{a}{q} \right| \leq n \frac{y^{-1} s}{q T^{\frac{1}{3}}},
\]
we have
\begin{align*}
\frac{s}{T^{\frac{1}{3}}} \sum_{sn \leq T^{\frac{1}{3}}} F(nsy) &= O\left(\frac{y^{-2}}{q}\right)+\frac{s}{T^{\frac{1}{3}}} \sum_{sn \leq T^{\frac{1}{3}}} F\left(n \frac{a}{q}\right)\\
&=O\left(\frac{y^{-2}}{q}\right)+O\left(\frac{qs}{T^{\frac{1}{3}}} \right)+\frac{1}{q} \sum_{j=0}^{q-1} F\left(\frac{ja}{q}\right)
\end{align*}
by the periodicity of $F$.  As $a$ is coprime to $q$, it does not play a role in the last average and can be dropped.  Furthermore,  for any $t \leq \frac{1}{q}$,
\[
F\left(\frac{j}{q}\right)=O\left(\frac{y^{-1}}{q}\right)+F\left(\frac{j}{q}+t\right),
\]
so
\begin{align*}
\frac{1}{q} \sum_{j=0}^{q-1} F\left(\frac{j}{q}\right) &= O\left(\frac{y^{-1}}{q}\right)+\frac{1}{q} \sum_{j=0}^{q-1} \int_0^1 F\left(\frac{j+t}{q} \right) dt \\
&= O\left(\frac{y^{-1}}{q}\right)+\int_0^1 F(t) dt.
\end{align*}
As both $y^{-2}q^{-1}$ and $qsT^{-\frac{1}{3}}$ are $O(y)$, this implies the claim.
\end{proof}
By Str\"ombergsson's result \cite{strombergsson}, 
\[
\left| y \int_0^{y^{-1}} f(\xi h(t)) dt - \int f \; d\mu_X \right| \ll y^\beta \ll (\delta^{-2} r^{-1})^\beta,
\]
so we see from Claim \ref{claim::qbig} that 
\[
\left| \frac{(1+\gamma) M^\gamma}{T^{\frac{1}{3}}}\sum_{n \leq T^{\frac{1}{3}}} f(\xi h((1+\gamma)M^\gamma n)) - \int f d\mu_X \right| \ll (\delta^{2} r)^{-\beta}
\]
unless there is a $q \leq y^{-3} \leq  r^{3} $ and $a$ coprime to $q$ such that 
\[
\left| (1+\gamma)M^\gamma y - \frac{a}{q} \right| \ll M^{\gamma}y^{-1}T^{-\frac{1}{3}} \leq  r T^{-\frac{1}{3}+\gamma} .
\]
To conclude the proof of Theorem \ref{powerequidist}, we just have to show that this is a very exceptional occurrence. 

Fortunately,  this is what one would expect: If we let 
\[
I_{q, a}:=\left\{v \in \R: \left| v - \frac{a}{q} \right| \leq r T^{-\frac{1}{3}+\gamma} \right\}
\]
denote the problem intervals for $q \leq r^3$ and $(a, q)=1$,  we note that they are proportional to $r T^{-\frac{1}{3}+\gamma}$.  Moreover, given distinct intervals $ I_{q_1, a_1}, I_{q_2, a_2}$,  the gap between them is at least of order $r^{-6}$,  as
\[
\left|\frac{a_1}{q_1}-\frac{a_2}{q_2}\right|\geq \frac{1}{q_1 q_2} \geq r^{-6}.
\]
As $r \ll T^{4\eps}$,  this means that the set $E:=\bigcup_{q \leq r^3, (a, q)=1} I_{q, a}$ makes up only a tiny proportion of the entire range.  Unless the function
\[
G(t):= \frac{t^\gamma R}{(t^{1+\gamma}-W)^2+1}=t^\gamma y
\]
is highly concentrated on a small part of its range,  our problem case $\left\{t \leq T: (1+\gamma) G(t) \in E \right\}$ will thus only occur on a negligible proportion of $[0,T]$. The claim below shows that $G$ does not behave in this unusual manner.
\begin{claim} \label{claim::dernormal}
For all but a $O(\delta + \delta^{-5} r^7 T^{-\frac{1}{3}+\gamma})$ proportion of $t \leq T$,  there does not exist $q \leq r^3$ such that 
\[
\left|(1+\gamma) G(t)-\frac{a}{q}\right| \leq  r T^{-\frac{1}{3}+\gamma}.
\]
\end{claim}
Before we show the claim,  we show how it implies Proposition \ref{prop::intcase}.  The claim implies that at most a small proportion of the intervals we split $[0,T]$ into when applying Taylor approximation will be bad; for the others, we know equidistribution from Claim \ref{claim::qbig}.  Collecting all the different error terms together,  
\[
\left|\frac{1}{T} \sum_{n \leq T} f(ph(n^{1+\gamma})) - \int f \; d\mu_X \right| \ll \delta + \delta^{-5} r^7 T^{-\frac{1}{3}+\gamma} + (\delta^{-2} r^{-1})^\beta,
\]
where the error terms come from, in that order,  Lemma \ref{lem::orbitapprox} and Claim \ref{claim::dernormal},  the contribution of the problem intervals $I_{q, a}$ on which the sequence $\xi (1+\gamma) M^\gamma n $ does not equidistribute in the periodic horocycle,  and the comparison with $\int f \; d\mu_X$ on the good intervals.  Setting $\delta=r^{-\frac{1}{10}}$ takes care of the first and third term,  while,  recalling that $r \ll T^{4\eps}$,  we can control the second term by setting $\varepsilon=\frac{1}{100}$.  This concludes the proof of Proposition \ref{prop::intcase} (and thus also the proof of Theorem \ref{powerequidist}) with only Claim \ref{claim::dernormal} left to be shown.
\end{proof}
\begin{proof}[Proof of Claim \ref{claim::dernormal}]
To show this claim, we use the following simple lemma, whose proof is left to the reader as an exercise.
\begin{lemma} \label{lem::derhelp}
Let $I \subset \R$ be an open interval and let $G \colon I \to \R$ be continuously differentiable such that $0<c \leq |G^\prime(t)| \leq C$ for all $t \in I$.  Let $\theta>0$ and let $a_1<b_1<a_2<\dots<a_{n-1}<b_{n-1}<a_n$ be real numbers with the property that $b_i-a_i \leq \theta (a_{i+1}-b_i)$ for all $1 \leq i<n$.  Then for $E:=(a_1, b_1) \cup \dots \cup (a_{n-1}, b_{n-1})$,
\[
|\{t \in I: G(t) \in E \}| \leq 2 \theta C c^{-1} |I|
\]
provided that $|I| \geq \theta C c^{-1}$.
\end{lemma}
To apply this to the function 
\[
G(t)= \frac{t^\gamma R}{(t^{1+\gamma}-W)^2+1}=t^\gamma y
\]
we are interested in,  we need to calculate its derivative.  We see that
\[
\frac{dy}{dt}(t)=-\frac{2(1+\gamma) t^\gamma (t^{1+\gamma} - W)R}{((t^{1+\gamma}-W)^2+1)^2}=-\frac{2y(1+\gamma) t^\gamma (t^{1+\gamma} - W)}{(t^{1+\gamma}-W)^2+1}
\]
and thus
\begin{align*}
G^\prime(t)= y t^{\gamma-1}  \left(\gamma - \frac{2(1+\gamma) t^{1+\gamma} (t^{1+\gamma} -W) }{(t^{1+\gamma}-W)^2+1} \right).
\end{align*}
We recall that in Lemma \ref{lem::orbitapprox} we exclude an interval $J_0$ of proportion $\delta$ to assure $r^{-1} \ll y \ll \delta^{-2} r^{-1}$. We also exclude a set $J_1$ comprised of two intervals of proportion $\delta$ to assure $t \geq \delta T$ and $\left|\frac{W}{t^{1+\gamma}}-1\right| \geq \delta$. This assures that $r^{-1} T^{\gamma-1} \ll y t^{\gamma-1} \ll \delta^{-3} r^{-1} T^{\gamma-1}$ on the range $[0,T] \backslash (J_0 \cup J_1)$.  If we can bound the expression in the bracket in a similar manner up to factors of powers of $\delta^{-1}$,  the claim will follow from Lemma \ref{lem::derhelp}.

To do this,  we note that for $t \in [0,T] \backslash J_1$,  
\[
\left| \frac{1}{(t^{1+\gamma}-W)^2+1}-\frac{1}{(t^{1+\gamma}-W)^2} \right|=O(\delta^{-4} T^{-4(1+\gamma)}),
\]
which implies 
\[
G^\prime(t)=y t^{\gamma-1} \left(\gamma + \frac{2(1+\gamma)}{\frac{W}{t^{1+\gamma}}-1} + O(\delta^{-4} T^{-2}) \right).
\]
We set $J_2:=\left\{t: \left|\frac{W}{t^{1+\gamma}}-\left(1-\frac{(2+\gamma)}{\gamma}\right)\right| \geq \delta \right\}$, which is the interval of proportion $\delta$  on which the second term roughly cancels out the first.  We then have that
\[
\delta \ll \gamma \left|\frac{W}{t^{1+\gamma}}-1\right|^{-1}  \left|\frac{W}{t^{1+\gamma}}-1+\frac{(2+\gamma)}{\gamma}\right|=\left|\gamma + \frac{2(1+\gamma)}{\frac{W}{t^{1+\gamma}}-1} \right| \ll \delta^{-1}
\]
on $[0,T] \backslash (J_1 \cup J_2)$, which implies that
\[
\delta r^{-1} T^{\gamma-1} \ll |G^\prime(t)| \ll \delta^{-4} r^{-1} T^{\gamma-1}
\]
on $[0,T] \backslash (J_0 \cup J_1 \cup J_2)$. We can now apply Lemma \ref{lem::derhelp} to each of the intervals left.  Recalling that each problem interval $I_{q, a}$ is of length $r T^{-\frac{1}{3}+\gamma}$ and the gap between any two successive intervals is of size at least $0.9 r^{-6}$, we find that
\[
\frac{1}{T}\left| \left\{t \in [0,T] \backslash (J_0 \cup J_1 \cup J_2): (1+\gamma)G(t) \in E \right\}\right|  \ll \delta^{-5} r^7 T^{-\frac{1}{3}+\gamma}
\]
where as before $E=\bigcup_{q \leq r^3, (a, q)=1} I_{q, a}$.
This shows Claim \ref{claim::dernormal}, which was the last missing piece in the proof of Theorem \ref{powerequidist}.

\end{proof}

\printbibliography

\end{document}